\documentclass[a4paper]{article}
\usepackage{amsmath}
\usepackage{amsbsy}
\usepackage{amsfonts}
\usepackage{amstext}
\usepackage{amssymb}
\usepackage{mathbbol}
\usepackage{wasysym}
\usepackage{mathtools}
\usepackage[amsmath,thmmarks]{ntheorem}		
\usepackage{stackrel}

\usepackage{tikz}
   \usetikzlibrary{matrix}
   \usetikzlibrary{calc}
   \usetikzlibrary{decorations,snakes}

\usepackage{float}
\usepackage{subfigure}

\usepackage{enumitem}
\usepackage{mathtools}

\usepackage[a4paper, left=2.5cm, right=2.5cm, top=2.8cm, bottom=3.0cm, bindingoffset=0cm]{geometry}

\usepackage{color}

\definecolor{grey}{rgb}{0.5,0.5,0.5}

\newcommand{\RR}{\ensuremath{{\mathbb R}}}
\newcommand{\QQ}{\ensuremath{{\mathbb Q}}}
\newcommand{\ZZ}{\ensuremath{{\mathbb Z}}}

\newcommand{\Mod}{\operatorname{mod}}

\newcommand{\calA}{\ensuremath{{\cal A}}}
\newcommand{\calO}{\ensuremath{{\cal O}}}

\newcommand{\kato}{\ensuremath{\calO}}

\renewcommand{\ker}{\ensuremath{\text{ker}}}

\newcommand{\spann}[1]{\ensuremath{\text{span}_{#1} \:}}
\renewcommand{\dim}{\ensuremath{\text{dim}}}

\newcommand{\grg}[1]{\ensuremath{\text{Grk}_{#1}}}

\newcommand{\maxid}{\ensuremath{{\mathfrak m}}}
\newcommand{\spec}{\ensuremath{\text{Spec}}}

\newcommand{\Sym}{\ensuremath{{\text{Sym}}}}

\newcommand{\anni}{\ensuremath{\text{Ann}}}
\newcommand{\supp}{\ensuremath{\text{Supp}\:}}

\newcommand{\wirk}{\ensuremath{{}_{^{_{\bullet}}}}}
\newcommand{\pair}[1]{\ensuremath{\langle #1 \rangle}}

\newcommand{\ua}{\ensuremath{^{\ast}}}

\renewcommand{\i}{\ensuremath{^{-1}}}

\newcommand{\p}{\ensuremath{^{\prime}}}

\newcommand{\ol}[1]{\ensuremath{\overline{#1}}}

\newcommand{\rar}{\rightarrow}
\newcommand{\irar}{\hookrightarrow}
\newcommand{\srar}{\twoheadrightarrow}

\newcommand{\gdw}{\ensuremath{\ \Leftrightarrow\ }}

\newcommand{\LieG}{\ensuremath{{\mathfrak g}}}
\newcommand{\Cartan}{\ensuremath{{\mathfrak h}}}

\newcommand{\Torus}{\ensuremath{{\mathfrak t}}}
\newcommand{\UE}{\ensuremath{{\mathfrak U}}}

\newcommand{\liesln}{\ensuremath{{\mathfrak sl}}_n}

\newcommand{\bchi}{\ensuremath{{B^{\chi}}}}
\newcommand{\bxchi}{\ensuremath{{B^{x\chi}}}}
\newcommand{\bchip}{\ensuremath{{B^{\chi\p}}}}
\newcommand{\subtor}{\ensuremath{\mathfrak g}}
\newcommand{\subtorh}{\ensuremath{\mathfrak h}}
\newcommand{\weyl}{\ensuremath{\calA}}

\newcommand{\brac}[1]{\ensuremath{\langle #1 \rangle}}

\newcommand{\EHP}{\ensuremath{\text{EHP}}}

\theoremstyle{plain}
\theoremseparator  {.}
\theoremheaderfont {\bfseries}
\theorembodyfont   {\normalfont}
\theoremnumbering  {arabic}
\makeatletter
\newtheoremstyle{normal}%
{\item[\hskip\labelsep \theorem@headerfont ##1\ ##2\theorem@separator]\normalfont}%
{\item[\hskip\labelsep \theorem@headerfont ##1\ ##2]{\theorem@headerfont (##3)}\theorem@separator\ \normalfont}
\newtheoremstyle{nonumber}%
{\item[\theorem@headerfont\hskip\labelsep ##1\theorem@separator]\normalfont}%
{\item[\theorem@headerfont\hskip \labelsep ##1]{\theorem@headerfont (##3)}\theorem@separator\ \normalfont}
\makeatother

\theoremstyle{normal}
\newtheorem{thm}               {Theorem}

\newtheorem{lemma}      [thm]  {Lemma}

\newtheorem{cor}        [thm]  {Corollary}

\newtheorem{defi}       [thm]  {Definition}
\newtheorem{prop}       [thm]  {Proposition}
\newtheorem{bem}        [thm]  {Remark}
\newtheorem{bsp}        [thm]  {Example}

\theoremstyle{nonumber}
\theoremsymbol{$\smiley$}
\newtheorem{bew}{Proof}
\theoremsymbol{}

\theoremstyle{nonumber}
\newtheorem{thmnn}{Theorem}

\theoremstyle{nonumber}
\theoremsymbol{$\smiley$}
\newtheorem{proofsketch} {Sketch of Proof}
\theoremsymbol{}

\begin{document}

\title{Goldie rank of primitive quotients via lattice point enumeration}
\author{Joanna Meinel and Catharina Stroppel}
\date{}
\maketitle

\begin{abstract}
Let $k$ be an algebraically closed field of characteristic $0$. In \cite{mvdb} the authors classify primitive ideals for rings
of torus invariant differential operators. This classification applies  in particular to subquotients of localized extended Weyl
algebras $\mathcal{A}_{r,n-r}=k[x_1,\ldots,x_r,x_{r+1}^{\pm1}, \ldots, x_{n}^{\pm1},\partial_1,\ldots,\partial_n]$ where it can be made explicit in terms of convex geometry. We recall these result and then turn to the corresponding
primitive quotients and study their Goldie ranks.
We prove that the primitive quotients fall into finitely many families whose Goldie ranks are given by a common quasi-polynomial and then realize these quasi-polynomials as Ehrhart quasi-polynomials arising from convex geometry.
\end{abstract}

\section*{Introduction}

Let $k$ be an algebraically closed field of characteristic $0$. In \cite{mvdb} the authors classify primitive ideals for rings
of torus invariant differential operators. It applies in particular to the Weyl algebras
$k[x_1,\ldots,x_n,\partial_1,\ldots,\partial_n]$ and their generalizations described in \eqref{eq:Weyl} below.
In these cases the classification can be used for a description of the primitive ideals in terms of convex
geometry. We recall these results and then turn to the corresponding primitive quotients and study their Goldie ranks. Recall that
the Goldie rank of a (prime noetherian) ring $R$ is defined as
$$\grg{}(R)=\max\{k\ |\ I_1\oplus\ldots\oplus I_k\text{ is a direct sum of nontrivial left ideals in }R\},$$
and measures the size of the quotient ring in some (non-commutative) way.

Based on the convex geometry used to classify the primitive ideals, we prove that the primitive quotients fall into finitely many families such that their Goldie ranks are given by a common quasi-polynomial.
Let us formulate our result more precisely: On the ring of differential operators on $k^r\times (k\ua)^s$
\begin{eqnarray}
\label{eq:Weyl}
\weyl=D(k[x_1,\ldots,x_r,x_{r+1}^{\pm1},\ldots,x_n^{\pm1}])=k[x_1,\ldots,x_r,x_{r+1}^{\pm1},\ldots,x_n^{\pm1},\partial_1,\ldots,\partial_n]
\end{eqnarray}
one has a $\ZZ^n$-grading coming from the adjoint action of the subspace $\Torus=\spann{k}\{x_i\partial_i\ |\ 1\leq i\leq n\}$
of $\weyl$. This grading carries over to the subquotient $\bchi = \weyl^\subtor /(\subtor-\chi(\subtor))$ of the algebra $\weyl$
for any subspace $\subtor\subset\Torus$ and $\chi\in\subtor\ua$, where $\weyl^\subtor$ denotes the $\subtor$-invariants of
$\weyl$ and $(\subtor-\chi(\subtor))$ stands for the ideal generated by elements $t-\chi(t)$ for $t\in\subtor$. The classification result of \cite{mvdb} states that
\begin{align*}
		\left\{\text{Primitive ideals of }\bchi\right\}
		& \stackrel{1:1}{\longleftrightarrow}\left\{\text{Regions }\ol{\brac{\alpha}}\subset\Torus\ua\ |\ \alpha\in
V(\subtor-\chi(\subtor))\subset\Torus\ua \right\}\\
		J(\alpha) & \longleftrightarrow \ol{\brac{\alpha}}
\end{align*}
where a region $\ol{\brac{\alpha}}$ is defined to be the Zariski closure of the support of a certain simple $\Torus\ua$-graded
module $L(\alpha)$ with annihilator $J(\alpha)$. In particular, every primitive ideal is isomorphic to the annihilator
$J(\alpha)$ of one of the $L(\alpha)$.
The regions $\ol{\brac{\alpha}}$ are given in terms of lattice points (coming originally from the weight lattice $\ZZ^n$ of
$\weyl$) and parallel translates of hyperplanes. This description gives a beautiful interpretation of the Goldie rank, namely it equals the number of connected components of
$\ol{\brac{\alpha}}$, see Theorem \ref{thm:conncpt}.

We are interested in the  behaviour of Goldie ranks under dilation
$\alpha\mapsto x\alpha$ with a positive integral factor $x\in\ZZ_{>0}$ and want to interpret this number in terms of counting integral
points in an appropriate polytope. For this purpose, we have to assume three technical conditions given in
Section~\ref{ssec:ehrhartalpha}. We then construct the polytopes so that the number of points in the intersection with the standard
$\ZZ$-lattice equals the connected components of $\ol{\brac{\alpha}}$. This counting of points can be done using Ehrhart theory for all $x\in\ZZ_{>0}$ such that $x\alpha$ satisfies the above technical conditions. The Goldie ranks for families of primitive quotients $\bchi/J(\alpha)$, see Theorem \ref{thm:goldpol}, can then be expressed in terms of the Ehrhart polynomial.

\begin{thmnn}
Let $\alpha\in\Torus\ua$ and $x\in\ZZ_{>0}$. Assume conditions \eqref{rat}-\eqref{mj} hold for $\alpha$ and $x\alpha$. Then the Goldie rank of
the primitive quotient $B^{x\chi}/J(x\alpha)$ is quasi-polynomial in $x$ and given by
$$\grg{}(B^{x\chi}/J(x\alpha))\ =\ \EHP_{Q}(x\p),$$
where $\EHP_{Q}$ is the Ehrhart quasi-polynomial of an appropriate rational polytope $Q$
with respect to the standard lattice and $x\p\in\QQ$ is the dilation factor obtained from rescaling $x$.
\end{thmnn}

This result should remind of the classical Goldie rank theory for universal enveloping algebras $\UE(\LieG)$ of semisimple
complex Lie algebras developed by Joseph.
First of all, for $\UE(\LieG)$ there is Duflo's theorem \cite{Duflo} that any primitive ideal arises as the annihilator of some
simple highest weight modules $L(\lambda)$ of highest weight $\lambda\in\Cartan\ua$ for a Cartan subalgebra
$\Cartan\subset\LieG$, see also \cite[Corollary 7.4]{Jantzen}. In other words, the BGG category $\kato\subset\UE(\LieG)$-$\Mod$
suffices for the description of the primitive ideals. It has consequences for the description of the corresponding primitive
quotients: Joseph found families of primitive ideals (depending on the central character respectively highest weight) in which the Goldie ranks vary
polynomially, \cite[Corollary 5.12]{Josephgr1}, see also \cite[Theorem 12.6]{Jantzen}. These Goldie rank polynomials can in
general only be computed up to a constant factor \cite[Appendix 14A.3]{Jantzen}.
Notice that for $\LieG=\liesln$ there is some recent progress, \cite{premet} and \cite{brundan}, using finite $W$-algebra
techniques, generalizing the fact that for a finite dimensional simple $\UE(\liesln)$-module, its dimension equals the Goldie
rank of its primitive quotient, \cite{Joseph1}.

The classification of primitive ideals was established in \cite{mvdb} in fact for a more general class of torus invariant
differential operators than the algebras studied here and resembling strongly the classical Lie theory by working wit the full subcategory $\kato^{(1)}$ of $\bchi$-modules $M$ such that
    $\bigoplus\limits_{\alpha\in\Torus\ua}\Sym(\Torus)/\maxid_\alpha \srar M$
as left $\Sym(\Torus)$-modules, where $M$ becomes a $\Sym(\Torus)$-module via $\Sym(\Torus)\subset\weyl^\subtor\srar \bchi$,
and the $\maxid_\alpha$ are the maximal ideals of $\Sym(\Torus)$. The modules in $\kato^{(1)}$ are $\Torus\ua$-graded with
$ M_\alpha = \{x\in M\ |\ \maxid_\alpha x=0\}$. The annihilators $J(\alpha)$ of the simple quotients $L(\alpha)$ of analogues of Verma modules give then a complete list of primitive ideals of $\bchi$. Unfortunately, the
aforementioned translation into convex geometry is however not available in general.

The first three sections contain basic definitions and recall results from \cite{mvdb} with the goal of formulating Theorem
\ref{thm:conncpt} which will then be proved in the last section. Section \ref{sec:geo} describes the link to the hyperplane geometry.
The new results appear in the last section (Theorems \ref{thm:goldpol} and \ref{thm:goldpol1}).

\subsubsection*{Acknowledgments} We are grateful to Gwyn Bellamy and Volodymyr Mazorchuk and the referee for helpful comments on a preliminary version of this paper. The second author thanks Ken Brown, Toby Stafford and Michel Van den Bergh for fruitful discussions.

\section{Preliminaries: localized extended Weyl algebras}\label{sec:wdh}

This section provides a brief overview of some results that can be found in \cite{mvdb}. As base field fix an algebraically closed field $k$ with $\text{char}(k)=0$.
\begin{defi}
Let $n,r,s\in\ZZ_{\geq 0}$ and $n=r+s$. The \emph{localized extended Weyl algebra} $\weyl=\weyl_{r,s}$ is defined to be the
$k$-algebra
	generated by $x_1,\ldots x_r$, $x_{r+1}^{\pm 1},\ldots,x_n^{\pm 1}$ and
	$\partial_1,\ldots,\partial_n$, subject to the relations
	$[x_i, x_j] = 0$, $[\partial_i, \partial_j] = 0$, $[\partial_i, x_j] = \delta_{ij}$ and $x_i x_i\i = 1$.
	By abuse of notation, we often write
	$\weyl=k[x_1,\ldots x_r,x_{r+1}^{\pm 1},\ldots,x_n^{\pm 1},\partial_1,\ldots,\partial_n]$.
\end{defi}
In the case $r=n$ this definition gives the classical Weyl algebra of polynomial differential operators. We
will usually call $\weyl_{r,s}$ just ${\emph Weyl\; algebra}$. It shares many properties with the classical Weyl algebra:
\begin{prop}\label{prop:weyltecheig}
	\begin{enumerate}[label=\roman{*}), ref=(\ref{prop:weyltecheig}.\roman{*})]
	\item\label{prop:weyltecheig-1} The algebra $\weyl$ can be written as the tensor product of
		$n$ small Weyl algebras $\weyl_i$ where $\weyl_i=k[x_i,\partial_i]$ for $i\leq r$ and
		$\weyl_i=k[x_i^{\pm1},\partial_i]$ for $i>r$,
		that is, we have an isomorphism of algebras $\weyl\cong\weyl_1\otimes\ldots\otimes\weyl_n$.
	\item $\weyl$ is a left (and right) noetherian domain,
	furthermore $\weyl$ is simple (with respect to twosided ideals).
	\end{enumerate}
\end{prop}
The proofs work exactly as for the classical Weyl algebra, see e.g. \cite{coutinho}; details can be found in \cite{dipl}.

\subsection{The weight lattice of the Weyl algebra}

Consider the subspace $\Torus = \spann{k}\{\pi_1, \ldots, \pi_n \ | \ \pi_i = x_i \partial_i \ \in\weyl\}$ of the Weyl algebra $\weyl$ with the adjoint action of $\Torus$ on $\weyl$ given by $[t,a]=ta-at$ for $t\in\Torus$, $a\in\weyl$.
Denote the weight spaces of $\weyl$ with respect to this action by $\weyl_\alpha$ with $\alpha\in\Torus\ua$.

\begin{lemma}\label{lem:weightlatticeweyl}
Identify $\Torus\ua$ with $k^n$ via $\pi_i\ua\mapsto e_i$ where $\pi_i\ua(\pi_j)=\delta_{ij}$.
Then
$$\weyl = \bigoplus_{\alpha \in \ZZ^n} \weyl_\alpha.$$
\end{lemma}
More precisely, the weight spaces are of the form $\weyl_\alpha=\weyl_0 \cdot a_\alpha$ with $\weyl_0=\Sym(\Torus) =
k[\pi_1,\ldots,\pi_n]$ and $a_\alpha = \prod\limits_{i=1}^n x_i^{(\alpha_i)}$, where $x_i^{(\alpha_i)}=x_i^{\alpha_i}$, if $i>r$
or $i\leq r$ and $\alpha_i \geq 0$, and $x_i^{(\alpha_i)}=\partial_i^{-\alpha_i}$, if  $i\leq r$ and $\alpha_i < 0$,
e.g. in $\weyl=k[x_1,x_2,x_3^{\pm1},\partial_1,\partial_2,\partial_3]$, one has
$a_{(-4,5,-6)} = \partial_1^4 x_2^5 x_3^{-6}$.
\begin{proofsketch}
A case by case calculation gives that $\Torus$ acts on $\weyl_0 \cdot a_\alpha$ by $\alpha$. To see that every monomial in
$\weyl$ lies in $\bigoplus_{\alpha \in \ZZ^n} \weyl_0 \cdot a_\alpha$, it suffices to consider the smallest Weyl algebras
$\weyl_i$ from Proposition \ref{prop:weyltecheig} because monomials in $\weyl_i$ commute with monomials in $\weyl_j$ for
$i\neq j$. To see that a monomial in $\weyl_i$ is indeed in $\bigoplus_{\alpha \in \ZZ^n} k[\pi_i]x_i^{(\alpha_i)}$, one
reorders the factors of the monomial, at the cost of summands with a lower number of factors $\partial_i$, and uses induction.
\end{proofsketch}

\subsection{The central subquotient $\bchi$}
\begin{bem}
The triple $(\weyl,\Torus,\phi)$ with $\phi=\text{incl}: \Torus\irar\weyl$ fits into a more general picture developed in
\cite[Section 3]{mvdb}. By taking subalgebras and quotients of such a triple one can construct new ones, see \cite[Section
4]{mvdb}. Here we are only interested in a combined construction resulting in a subquotient of the Weyl algebra $\weyl$.
\end{bem}
\begin{defi}
	Let $\subtor\subset\Torus$ be a subspace and let $\chi \in \subtor\ua$. Then define
	$$ \bchi \ = \ \weyl^\subtor / (\subtor - \chi(\subtor))$$
	where $\weyl^\subtor = \{ a \in \weyl \ |\ \forall t \in \subtor: \ [t,a] =0 \}$
	is the centralizer for the action of $\subtor \subset \Torus$.
	The twosided ideal $(\subtor - \chi(\subtor))$ of $\weyl^\subtor$ is generated by
	\mbox{$\subtor - \chi (\subtor) = \{ t - \chi(t) \ \in \weyl \ |\ t \in \subtor\}$}.
\end{defi}
On $\bchi$, we still have an action of $\Torus$ coming from the precomposition of
$\phi:\ \Torus\subset \weyl^\subtor\srar\bchi$ with the adjoint action; $\bchi$ inherits then a weight space decomposition from $\weyl$. Note that in $\weyl^\subtor$ only the weight spaces for $\alpha\in V(\subtor)$ survive, i.e.
$\weyl^\subtor = \bigoplus_{\alpha \in V(\subtor)} \weyl_\alpha$.

\section{Primitive ideals and simple modules}

We turn now to the description of the primitive ideals and the primitive quotients of the algebra $\bchi$. The primitive ideals of $\bchi$ have a very nice classification given in \cite[Theorem 7.3.1]{mvdb} which reduces the study of
\emph{all} the primitive ideals to the study of annihilators of a very special class of simple modules.
\begin{thm}[Correspondence between primitive ideals and closed regions]\label{thm:einszueins}
There is a one-to-one correspondence
		\begin{align*}
		\left\{ \text{ Regions } \ol{\brac{\alpha}}\subset\Torus\ua\ | \ \alpha \in V(\ker \phi) \right\}
		& \stackrel{1:1}{\longleftrightarrow} \left\{ \text{ Primitive ideals } \subset \bchi\right\} \\
		\ol{\brac{\alpha}} &\longleftrightarrow J(\alpha)
		\end{align*}
where $J(\alpha)$ is the annihilator of the simple $\Torus\ua$-graded module $L(\alpha)$ whose construction is given in the
subsequent section, and $\brac{\alpha}$ is the support of $L(\alpha)$, i.e. the set of weights $\beta\in\Torus\ua$ for which $L(\alpha)_{(\beta)}\neq0$, with Zariski closure
$\ol{\brac{\alpha}}\subset\Torus\ua$. The map $\phi$ is given by $\Torus\subset \weyl^\subtor\srar\bchi$.
\end{thm}
\begin{bem}
This theorem holds for the more general class of rings studied in \cite[Theorem 3.2.4]{mvdb}. The
construction of the simple graded module $L(\alpha)$ (outlined below) works as in the general setup.
\end{bem}

\subsection{Some simple graded modules for $\bchi$}\label{ssec:lalpha}

Denote the maximal ideal corresponding to $\alpha\in\Torus\ua = \spec(\Sym(\Torus))$ by $\maxid_\alpha$. We consider
$\maxid_\alpha$ as subset of $\bchi$ via $\Torus\stackrel{\phi}{\rar}\bchi$. By \cite[Proposition 3.1.7]{mvdb}, the left $\bchi$-module
$\bchi/\bchi\maxid_\alpha$ comes with a $\Torus\ua$-grading:
\begin{lemma}
The $\bchi$-module $\bchi/\bchi\maxid_\alpha$ is graded with
\begin{align*}
(\bchi/\bchi\maxid_\alpha)_{(\gamma)}\ :=\ \{m\in \bchi/\bchi\maxid_\alpha\ |\ (\maxid_\gamma)\wirk m =0\}\ =\
B^\chi_{\gamma-\alpha}/B^\chi_{\gamma-\alpha}\maxid_\alpha.
\end{align*}
Furthermore, $\ol{1}\in \left(\bchi/\bchi\maxid_\alpha\right)_{(\alpha)}$, so $\bchi/\bchi\maxid_\alpha$ is generated in degree
$(\alpha)$.
\end{lemma}
\begin{bem}
One can check that all weights of $\bchi/\bchi\maxid_\alpha$ are contained in $V(\ker(\phi))\subset\Torus\ua$.
\end{bem}
This grading plays the role of a weight space decomposition in the analogy with $\UE(\LieG)$, in which context the module
$\bchi/\bchi\maxid_\alpha$ appears as `highest weight module of weight $\alpha$'.
\begin{prop}[The unique simple quotient $L(\alpha)$ of $\bchi/\bchi\maxid_\alpha$]\label{prop:317einf}
\begin{enumerate}[label=\roman{*}), ref=(\ref{prop:317einf}.\roman{*})]
	\item\label{prop:317einf1}
		$\bchi/\bchi\maxid_\alpha$ has a unique simple quotient $=: \ L(\alpha)$
	\item\label{prop:317einf3}
		$\dim_{k} \left(\bchi/\bchi\maxid_\alpha\right)_{(\beta)} \leq 1$
		and therefore $\dim_{k} L(\alpha)_{(\beta)}\leq 1$ for all $\beta\in V(\ker\phi)$.
	\item\label{prop:317einf5}
		$L(\alpha_1) \cong L(\alpha_2)$ if and only if $\supp L(\alpha_1) \cap \supp L(\alpha_2) \neq \emptyset$.
	\end{enumerate}
\end{prop}
A very important feature of the module $L(\alpha)$ is that it inherits the $\Torus\ua$-grading from $\bchi/\bchi\maxid_\alpha$.
\begin{bem}
Property \ref{prop:317einf5} means that knowing one weight of the special simple module determines it already up to isomorphism.
In other words, if $\beta\in L(\alpha)$, then $L(\beta)\cong L(\alpha)$.
\end{bem}
\begin{defi}
Let $\alpha\in V(\ker(\phi))$. Define the region
$$\brac{\alpha}_\bchi\ =\ \supp L(\alpha)\ =\ \{\beta\in V(\ker(\phi))\ |\ L(\beta)\cong L(\alpha)\}.$$
and the primitive ideal
$$J(\alpha)\ =\ \anni(L(\alpha))\ =\ \{a\in A\ |\ a\wirk L(\alpha)=0\}.$$
\end{defi}
\begin{bem}\label{rem:vkerphibchi}
For our $\phi$ given by $\Torus\subset \weyl^\subtor\srar\bchi$, one can check that
$V(\ker(\phi))= V(\subtor-\chi(\subtor))$.
\end{bem}

\subsection{The regions of the Weyl algebra $\weyl$ and its subquotient $\bchi$}\label{ssec:reg}

We want to give a detailed description of the regions $\brac{\alpha}_\bchi$ and their Zariski closures deduced from a description of the regions $\brac{\alpha}_\weyl$, which are by definition the supports of the simple modules $L(\alpha)$.

One has $\brac{\alpha}_\weyl\subset \alpha + \{\text{ weights of }\weyl\}$ for any $\alpha\in\Torus\ua$, because $L(\alpha)$ is
a graded $\weyl$-module generated in degree $\alpha$. In Lemma \ref{lem:weightlatticeweyl} it turned out that the weights of
$\weyl$ form a lattice $\ZZ^n$ inside $\Torus\ua$, in other words, $\brac{\alpha}_\weyl\subset \alpha +\ZZ^n$. This plays indeed
an important role in the description of the regions $\brac{\alpha}_\weyl$, as the following result shows. It is a slight
reformulation of \cite[Corollary 6.2]{mvdb}.

\begin{prop}\label{prop:bracalpha}
For all $\alpha\in\Torus\ua$, the region $\brac{\alpha}_\weyl$ is of the form
$$\brac{\alpha}_\weyl = \{ \beta\in \Torus\ua\ |\ \beta\equiv\alpha\mod\ZZ^n\ \text{ and }
	\beta_i\in\ZZ_{\geq0} \gdw \alpha_i \in \ZZ_{\geq0}\ \forall 1\leq i\leq r \}.$$
\end{prop}
The regions behave well under the passage from $\weyl$ to $\bchi$. Indeed, by \cite[Proposition 4.4.1]{mvdb},
\begin{eqnarray}
\label{eq:bracalphabchi}
\brac{\alpha}_\bchi &=& \brac{\alpha}_\weyl \cap V(\subtor-\chi(\subtor))
\end{eqnarray}
holds for $\alpha \in V(\subtor-\chi(\subtor))$. We have already a very concrete description of $\brac{\alpha}_\weyl$, here is
one for $V(\subtor-\chi(\subtor))$:
\begin{lemma}\label{lem:vg-chig}
Define $\eta_i = \pi_i\ua|_\subtor$.
Under the usual identification of $\Torus\ua$ with $ k^n$,
$$V(\subtor-\chi(\subtor))\ =\{\alpha=(\alpha_1,\ldots \alpha_n)\in k^n\ |\ \sum_{i=1}^n\alpha_i\eta_i =\chi\ \},$$
thus $V(\subtor-\chi(\subtor))=V(\subtor)+\alpha$ for any $\alpha\in V(\subtor-\chi(\subtor))$.
\end{lemma}
Putting together the results of Proposition \ref{prop:bracalpha} and \eqref{eq:bracalphabchi}, one can see how the geometrical
picture evolves: $\brac{\alpha}_\bchi$ is the intersection of $V(\subtor-\chi(\subtor))$ with the lattice $\alpha+\ZZ^n$ and the
polyhedral cone
$$\Delta_\alpha = \{\beta\ |\ \beta_i\in\ZZ_{\geq0} \gdw \alpha_i \in \ZZ_{\geq0}\text{ for }1\leq i\leq r\}.$$
\begin{bem}
Note that from Proposition \ref{prop:bracalpha} it follows that the localization of $x_i$ leads to \emph{less} inequalities in
the description of the regions. E.g. for $r=0$ where all $x_i$ are inverted, no inequality is left and
$\brac{\alpha}_\weyl=\alpha+\ZZ^n$. In other words, the concrete geometrical classification of the primitive ideals for any
subquotient $\bchi$ with $r=0$ tells us that there is only one such ideal, corresponding to
$V(\subtor-\chi(\subtor))\subset\Torus\ua$, and this is the zero ideal. So $\bchi$ is simple in this case, see \cite[Proposition
3.3.1]{mvdb}.
\end{bem}

\section{Primitive Quotients and Goldie Rank}

The goal of this section is a description for the Goldie rank of the primitive quotients of $\bchi$ involving the geometry of
the regions. We recall the definition of Goldie rank and the description of the primitive quotients of $\bchi$ and formulate the main theorem which will be proved finally in Section \ref{sec:goldie}.

We work with the following definition of Goldie rank:
\begin{defi}
	For a prime noetherian ring $R$, the Goldie rank is defined to be
	$$\grg{}(R)\
	=\ \max\{n\ |\ \text{there is a direct sum of left ideals }I_1\oplus\ldots\oplus I_n\subset R\}.$$
\end{defi}
\begin{bem}
This definition applies since the ring $\bchi/J$ constructed from the Weyl algebra is indeed prime and noetherian: As a primitive
ring, $\bchi/J$ is also prime \cite[Proposition 11.6]{LamFC}. Moreover, $\bchi$ is noetherian (see \cite[Theorem 7.3.1]{mvdb})
and so is $\bchi/J$.
\end{bem}
Note that if the prime noetherian ring $R$ has no zero divisors, then $\grg{}(R)=1$. This follows from the alternative characterization of Goldie rank, given as the size $n$ of the matrix ring over the skew field $D$ such that the classical ring of quotients $Q$ of $R$ is isomorphic to $M^n(D)$, but in this case $Q=D=M^1(D)$.

\subsection{Goldie rank as number of connected components}

A {\emph primitive quotient} of a ring $R$ is defined to be the quotient of $R$ by a primitive ideal $J$. From Theorem
\ref{thm:einszueins} we have a complete list of the primitive ideals in $\bchi$, they are all of the form $J(\alpha)$ for
$\alpha\in V(\ker(\phi))$. \cite[Proposition 7.4.1]{mvdb} provides a description of the corresponding primitive quotient
$\bchi/J(\alpha)$, which we recall now. Because one has to distinguish between the $\bchi$ with underlying space $\subtor$ and
some other $\bchip$ corresponding to $\subtorh \supset\subtor$, we index it by the underlying space.

\begin{prop}[\cite{mvdb}]
\label{prop:primquot}
	Let $J\subset B^\chi_\subtor$ be a primitive ideal.
	Then there is a subspace $\subtor\subset\subtorh\subset\Torus$ and
	$\chi_1,\ldots,\chi_p\in \subtorh\ua$ with $\chi_i |_\subtor =\chi$ for all $i$ such that
	$$B^\chi_\subtor/J \ = \
	\left(
	\begin{array}{cccc}
		B_\subtorh^{\chi_1} &B_\subtorh^{\chi_1,\chi_2} &&\\
		B_\subtorh^{\chi_2,\chi_1} &B_\subtorh^{\chi_2}&&\\
		&&\ddots&\\
		&&&B_\subtorh^{\chi_p}
	\end{array}
	\right),
$$
where $B_\subtorh^{\chi_i,\chi_j} = \left(B^\chi_\subtor\right)^{\chi_i-\chi_j}_\subtorh /
(\subtorh-\chi_i(\subtorh))\left(B^\chi_\subtor\right)^{\chi_i-\chi_j}_\subtorh $ is a bimodule with
$\left(B^\chi_\subtor\right)^{\chi_i-\chi_j}_\subtorh\
=\ \{ a\in B^\chi_\subtor\ |\ [\phi(t),a] = (\chi_i-\chi_j)(t)a\ \forall \ t\in\subtorh\}.$
The right hand side consists of matrices with entries $b_{ij}\in B^{\chi_i,\chi_j}$, the multiplication is given by `matrix
multiplication'.
\end{prop}

The following theorem connects the Goldie rank of the primitive quotients with the geometry of the closed regions
$\ol{\brac{\alpha}}_\bchi$, crucial for the upcoming theorem on Goldie rank quasi-polynomials. The statement can be
found (without proof) in \cite[Corollary 7.4.3]{mvdb} and will be proved in Section ~\ref{sec:goldie}.

\begin{thm}
\label{thm:conncpt}
	The Goldie rank of $\bchi/J(\alpha)$
	equals the number of connected components of $\ol{\brac{\alpha}}_\bchi$.
\end{thm}

\section{Primitive quotients of $\bchi$ and hyperplane arrangements}\label{sec:geo}

The goal of this section is to give a geometrical description of the Zariski closed regions $\ol{\brac{\alpha}}_{\bchi}\subset\
k^n$ in terms of hyperplane arrangements.

\subsection{Computation of the closure of $\brac{\alpha}_\bchi$}\label{ssec:closurereg}

We have already a description of the non-closed regions $\brac{\alpha}_\bchi$ as the intersection of $V(\subtor-\chi(\subtor))$ with
the lattice $\alpha+\ZZ^n$ and the polyhedral cone $\Delta_\alpha$. The cone cuts out the coordinate chamber in which $\alpha$
lives, but only for the `interesting' indices $\{ i\ |\ 1\leq i\leq r,\ \alpha_i\in \ZZ\} =: T_\alpha$. Indices that are not in
$T_\alpha$ are not affected.

To compute the Zariski closure of $\brac{\alpha}_\bchi$, one may apply the following result from convex geometry (see
\cite[Proposition 7.1.2]{mvdb}). For this purpose, we assume from now on that the arrangement is \emph{rational},
i.e. the defining equations of $V(\subtor-\chi(\subtor))\subset k^n$ have coefficients in $\QQ^n$ and in particular
$\alpha\in\QQ^n\subset k^n$.

\begin{prop}\label{prop:geo}
	Let $E$ be a vector space over $\QQ$.
	Let $L\subset E$ be a full $\ZZ$-lattice, i.e. $L$ generates $E$.
	\begin{enumerate}[label=\roman{*}), ref=(\ref{prop:geo}.\roman{*})]
	\item\label{lem:techgeo} Given any $\lambda_1,\ldots,\lambda_m \ \in \ E\ua$,
	there is a unique decomposition of the index set $T=\{1,\ldots,m\}$ into two disjoint parts $I\dot{\cup}J$, such that
	there are $e \in E, \ z = (z_1,\ldots,z_m) \in \QQ^m$ with
\begin{equation*}
	\sum\limits_{i=1}^m z_i \lambda_i \ = \ 0, \;\text{   and   }\;\langle \lambda_i, e\rangle = \lambda_i(e)= \
	\begin{cases} >0, &\text{for } i \in I\\ =0, &\text{for } i \in J \end{cases}\;\text{   and   }\;
	z_i = \ \begin{cases} =0, &\text{for } i \in I\\ >0, &\text{for } i \in J. \end{cases}	
    \end{equation*}
	\item Given furthermore $q_1,\ldots,q_m\in\QQ$, define $E\p=\bigcap_{j\in J} \ker(\lambda_j)$ and
	\begin{eqnarray*}
	C=\{x\in E \;|\; \langle \lambda_i,x\rangle=\lambda_i(x) \leq q_i, \; \forall i\in T\},
	&&C\p=\{x\in E \; |\; \langle \lambda_j,x\rangle=\lambda_j(x) \leq q_j \; \forall j\in J\},
	\end{eqnarray*}
	then the Zariski closure of $C\cap L$ equals
	$\ol{C\cap L} = C\p \cap (L+E\p)$
	and
\label{prop:geo-2}
	$C\p \cap(L+E\p)$ is a finite union of translates of $E\p$.
	\end{enumerate}
\end{prop}

To apply the proposition to $\brac{\alpha}_\bchi = V(\subtor-\chi(\subtor))\cap (\alpha+\ZZ^n)\cap\Delta_\alpha$, we have to
consider the lattice and the cone $\Delta_\alpha = \{\beta\ |\ \beta_i\in\ZZ_{\geq0} \gdw \alpha_i \in \ZZ_{\geq0}\text{ for
}1\leq i\leq r\}$ over $\QQ$. Moreover, everything has to be translated by $-\alpha$ back to the origin, i.e. the proposition
will give us a description of the closure of $\QQ^n\cap\brac{\alpha}_\bchi -\alpha$.
Define
\begin{equation*}
E\ =\ V(\subtor)\cap\QQ^n,\quad L\ =\ \ZZ^n\cap V(\subtor),\quad T\ =\ T_\alpha \ =\ \{ i\ |\ 1\leq i\leq r,\ \alpha_i\in \ZZ\},\\
\end{equation*}
\begin{eqnarray*}
\lambda_i\ =\
	\begin{cases}\lambda_i(u)= -u_i\ &\text{if }\alpha_i\in\ZZ_{\geq0},\\
	\lambda_i(u)=u_i\ &\text{if }\alpha_i\in\ZZ_{<0}\end{cases},
&&
q_i\ =\
	\begin{cases} \alpha_i\ &\text{if }\alpha_i\in\ZZ_{\geq0},\\
	-\alpha_i-1\ &\text{if }\alpha_i\in\ZZ_{<0},\end{cases}
\end{eqnarray*}
so that $C\ =\ \{x\in E \ |\ \lambda_i(x) \leq q_i \text{ for all } i\in T\}=\QQ^n\cap\Delta_\alpha -\alpha$, and thus we have
indeed
$$\QQ^n\cap\brac{\alpha}_\bchi -\alpha\ =\ L\cap C.$$
Now Proposition \ref{prop:geo} tells us that $\ol{C\cap L}\ =\ C\p \cap (L+E\p)$ with
\begin{align*}
C\p&:=\ \{x\in E \ |\ \lambda_j(x) \leq q_j, \text{ for all }j\in J\}=\ \{\gamma\in V(\subtor)\ |\ \gamma_j\geq-\alpha_j \gdw
\alpha_j \geq0\text{ for all }j\in J\}\ \cap\ \QQ^n,\\
E\p&:=\ \bigcap_{j\in J} \ker(\lambda_j)=
\ \{\gamma\in V(\subtor)\ |\ \gamma_j=0, \text{ for all } j\in J\subset T_\alpha\}\
\cap\QQ^n.
\end{align*}
Here, the index set $J$ is chosen by Proposition \ref{lem:techgeo} in such a way that the corresponding inequalities cut out a
\emph{bounded} polyhedron. The other indices either belong to $T_\alpha\setminus J$, so they give inequalities that are no
longer visible after the closure, or they do not even belong to $T_\alpha$, in which case there has never been associated any inequality to them. 
From now on we denote by $I$ all indices that are not in $J$, $I=\{1,\ldots,n\}\setminus J$, not only the index set
$T_\alpha\setminus J$ from Proposition \ref{lem:techgeo}.
\begin{cor}[Description of $\ol{\brac{\alpha}}_\bchi$]
Let $\alpha\in\Torus\ua$. The closure of $\brac{\alpha}_\bchi=\supp L(\alpha)$ is given by
	\begin{align*}
	\ol{\brac{\alpha}}_{\bchi}\
	&=\ \left\{ \gamma\in V(\subtor-\chi(\subtor))\ \left|\
		\begin{matrix}
			\text{for all }j\in J\text{ one has }\gamma_j\in\ZZ,\\
			\gamma_j\geq0\gdw\alpha_j \geq0,\\
			\gamma_j<0\gdw\alpha_j <0
		\end{matrix}\right.\right\}\\
	&\cap\
	\left\{ \gamma\in V(\subtor-\chi(\subtor))\ \left|\
		\begin{matrix}
			\sum\limits_{j\in J}\gamma_j\eta_j\ \in\ \sum\limits_{j\in J}\alpha_j\eta_j\ +\ \sum\limits_{i\in I}\ZZ\eta_i
	\end{matrix}\right.\right\}.
	\end{align*}
\end{cor}

\begin{bew}
Apply the proposition to the rational arrangement, then pass to the closure over $k$ and translate the result by $\alpha$ back
to $V(\subtor-\chi(\subtor))$. This gives
\begin{align*}
\ol{\brac{\alpha}}_{\bchi}\
&=\ \left\{ \gamma\in V(\subtor-\chi(\subtor))\ \left|\
	\begin{matrix}
		\exists \delta\in V(\subtor-\chi(\subtor)):\quad \delta\equiv\alpha\mod\ZZ^n,\\
		\text{and for all }j\in J,\\
		\gamma_j=\delta_j\geq0\gdw\alpha_j \geq0,\\
		\gamma_j=\delta_j<0\gdw\alpha_j <0
	\end{matrix}
	\right.\right\}
\end{align*}
where one has to notice that the inequalities still make sense because they only involve expressions in $\QQ$ (due to the
rationality assumption). The rest of the proof is calculation.
\end{bew}
\begin{bem}\label{bem:mjmtheta}
Interestingly, the first set in this intersection
$$\small
\left\{ \gamma\in V(\subtor-\chi(\subtor))\ \left|\
		\begin{matrix}
			\text{for all }j\in J\text{ one has }\gamma_j\in\ZZ,\\
			\gamma_j\geq0\gdw\alpha_j \geq0,\\
			\gamma_j<0\gdw\alpha_j <0
		\end{matrix}\right.\right\}
= \left\{ \gamma\in V(\subtor-\chi(\subtor))\ \left|\
		\begin{matrix}
			\text{for all }j\in J\text{ one has }\gamma_j\in\ZZ,\\
			\gamma_j\geq0\text{ for all }j\in J_+,\\
			\gamma_j<0\text{ for all }j\in J_-
		\end{matrix}\right.\right\} $$
does only depend on the sign configuration of $\alpha$, not on $\alpha$ itself. Thus, this part of the description is the same
for all $\alpha\in\Torus\ua$ with the same signs
\begin{align*}
J_+\ &=\ \{j\in J\ |\ \alpha_j\geq 0\}\\
J_-\ &=\ \{j\in J\ |\ \alpha_j <0\}.
\end{align*}
Different $\alpha$'s that share the same sign configuration $J=J_+\cup J_-$ do result in the same set, which we will denote by
$M_J$. This is in contrast to the case of the set
\begin{eqnarray}
\label{eqMtheta}
	\left\{ \gamma\in V(\subtor-\chi(\subtor))\ \left|\
		\begin{matrix}
			\sum\limits_{j\in J}\gamma_j\eta_j\ \in\ \sum\limits_{j\in J}\alpha_j\eta_j\ +\ \sum\limits_{i\in I}\ZZ\eta_i
	\end{matrix}\right.\right\},
\end{eqnarray}
which may change with $\alpha$, even for $\alpha$'s with the same sign configuration, \cite[Example 7.2.7]{mvdb}.
\end{bem}
Following \cite{mvdb}, we denote the equivalence class of $\sum\limits_{j\in J}\alpha_j\eta_j\ +\ \sum\limits_{i\in I}\ZZ\eta_i$
by $\vartheta$ and the set \eqref{eqMtheta} by $M_\vartheta$. With this shorthand, we have
\begin{eqnarray}
\label{eq:shorthand}
\ol{\brac{\alpha}}_\bchi&=&M_J\cap M_\vartheta.
\end{eqnarray}
Note that there are infinitely many different regions $\brac{\alpha}_\bchi$ (for arbitrary $\alpha\in\Torus\ua$, we even get
infinitely many different lattices $\alpha+\ZZ^n$ as long as $n>0$). However we recall the proof of the following result from
\cite[Corollary 7.2.6]{mvdb} to illustrate the geometry:
\begin{prop}\label{prop:zaehlt}
There are only finitely many regions closures $\ol{\brac{\alpha}}_\bchi$.
\end{prop}
\begin{bew}
There are at most $2^n$ many different sign configurations $J=J_+\cup J_-$ and so, thanks to \eqref{eq:shorthand}, we only have
to show that for a given sign configuration, the number of admissible cosets $\vartheta$ is also finite. Recall that $\vartheta
=\sum_{j\in J}\alpha_j\eta_j+\sum_{i\in I}\ZZ\eta_i$ belongs to a point $\alpha$ in $V(\subtor-\chi(\subtor))$ with sign
configuration in $J$. Therefore $\vartheta$ has to satisfy $\alpha_j\in\ZZ_{\geq0}$ for all $j\in J_+$ and $\alpha_j\in\ZZ_{<0}$
for all $j\in J_-$, furthermore there are $\mu_i\in\QQ$ such that $\sum_{j\in J}\alpha_j\eta_j\ +\ \sum_{i\i I}\mu_i\eta_i\ =\
\chi$.
We want to construct now $\psi\in\subtor$ and apply it to $\sum_{j\in J}\alpha_j\eta_j\ +\ \sum_{i\i I}\mu_i\eta_i\ =\ \chi$ to
show that there are only finitely many solutions $\alpha_j\in\ZZ$, hence only finitely many admissible $\vartheta$.
Use once more Proposition \ref{prop:geo} to obtain coefficients $z_k\in\QQ$ with
$$\sum_{k\in T_\alpha}z_k\lambda_k\ =\ 0\quad\text{and}\quad z_j>0\text{ for all }j\in J,\quad z_k=0\text{ otherwise.} $$
Now define $\psi\in\subtor=(\subtor\ua)\ua$ by $\brac{\psi,\eta_k}=y_k$ with $y_k=-z_k$ for $ k \in J_+$, $y_k=z_k$ for $k\in
J_-$ and $y_k=0$ else. This is well-defined: The $\eta_k$ generate $\subtor\ua$ but they need not be linearly independent, so we
have to check that $\sum_{k=1}^n \gamma_k y_k=0$ whenever $\sum\limits_{k=1}^n \gamma_k\eta_k=0$.
Now use that $\gamma=(\gamma_k)\in V(\subtor)$, so we may apply $0=\sum_{k\in T_\alpha}z_k\lambda_k \in V(\subtor)\ua$ to it and
find that
$$0=\sum_{k\in T_\alpha}z_k\lambda_k(\gamma)=\sum\limits_{k=1}^n \gamma_k y_k.$$
So $\psi$ exists and we find indeed the equation $\sum_{j\in J}\alpha_j y_j\ =\ \pair{\psi,\chi}$ with
$\alpha_j y_j\leq0$ and $\alpha_j\in\ZZ$ for all $j\in J$, which has only finitely many solutions.
\end{bew}

\begin{bsp}
In case $J=\{1,\ldots,n\}$, which occurs for the classical Weyl algebra
$k[x_1,\ldots,x_n,\partial_1,\ldots,\partial_n]$ only, one has $M_\vartheta=V(\subtor-\chi(\subtor))$, so $M_\vartheta$ doesn't
influence the region $\ol{\brac{\alpha}}_\bchi=M_J$. Furthermore, $E\p=\{0\}$, so in this case the region consists of finitely
many points. The case $J=\emptyset$ implies immediately $\ol{\brac{\alpha}}_\bchi= V(\subtor-\chi(\subtor))$.
\end{bsp}
This example leads to the question about inclusion relations between the $\ol{\brac{\alpha}}_\bchi$. In our notation, the answer
of \cite[Proposition 7.2.5]{mvdb} reads as follows:
\begin{prop}
Fix $\subtor\subset\Torus$. Let $\alpha,\beta$ be in $\Torus\ua$ with sign configurations $J^\alpha,J^\beta$ and corresponding
$\vartheta_\alpha,\vartheta_\beta$ from Remark \ref{bem:mjmtheta}.
Then
\begin{eqnarray*}
\ol{\brac{\alpha}}_\bchi\subseteq\ol{\brac{\beta}}_\bchi&\text{if and only if}&
\left(J_+^\beta\subseteq J_+^\alpha, \quad
J_-^\beta\subseteq J_-^\alpha,\quad\text{and}\quad
\vartheta_\beta=\vartheta_\alpha\mod\sum\limits_{i\in I^\beta}\ZZ\eta_i\right).
\end{eqnarray*}
\end{prop}

\subsection{The connected components of $\ol{\brac{\alpha}}_\bchi$}\label{ssec:zshgkomp}

From the application of Proposition \ref{prop:geo-2} to $\ol{\brac{\alpha}}_\bchi$ we also know a description of
$\ol{\brac{\alpha}}_\bchi$ as a finite union of translates of the hyperplane $E\p = \bigcap\limits_{j\in J} \ker(\lambda_j)$, in
other words,
$$\ol{\brac{\alpha}}_\bchi=\bigcup\limits_{i=1}^p (E\p+\delta_i +\alpha)$$
for some $\delta_i\in\Torus\ua$.
 For later purposes we examine this description in more detail. The following lemma is part of the proof of
 \cite[Proposition 7.4.1]{mvdb}. With $C$, $L$ etc. we use the same notation as before.
\begin{lemma}\label{lem:vsubtorh}
	Let $\alpha\in V(\subtor-\chi(\subtor))$.
	Then there are $\subtor\subset\subtorh\subset\Torus$ and $\chi_i\in\subtorh$ such that
	$$\ol{\brac{\alpha}}_\bchi=\ol{L\cap C\ +\ \alpha} =\bigcup_{i=1}^p E\p+(\delta_i+\alpha)
	= \bigcup_{i=1}^p V(\subtorh-\chi_i(\subtorh)),$$
	where $\subtorh= \spann{k}\{\lambda_j\ |\ j\in J\}$
	(under the canonical identification $(\Torus\ua)\ua\cong\Torus$
	such that $\lambda_j(\gamma)=\gamma(\lambda_j)$ for all $\gamma\in\Torus\ua$)
	and $\chi_i=(\delta_i+\alpha)|_\subtorh$.
\end{lemma}
We consider $V(\subtorh-\chi_i(\subtorh))$ and $V(\subtor-\chi(\subtor))$ as subsets of $\Torus\ua$ as usual.
\begin{bem}\label{bem:zshgkomp}
	From $\ol{\brac{\alpha}}_\bchi= \bigcup_{i=1}^p (\chi_i + V(\subtorh))$
	one can see that $\ol{\brac{\alpha}}_\bchi$ has exactly $p$ connected components:
	$V(\subtorh)$ is connected and the different translates have trivial pairwise intersection.
	Hence, the connected components are exactly the $V(\subtorh-\chi_i(\subtorh))$.
\end{bem}

\section{Goldie rank polynomials for primitive quotients of $\bchi$}\label{sec:goldie}

We finally use the results from Section \ref{sec:geo} to calculate the Goldie rank of the primitive quotient $\bchi/J(\alpha)$.
We establish the existence of Goldie rank polynomials for the central quotients of the Weyl algebras, and
connect them with the geometry of hyperplane arrangements studied in the previous section.

\subsection{Proof of Theorem \ref{thm:conncpt}}
Denote $B_\subtor^\chi/J(\alpha)$ by $B$. We have to consider left ideals in $B$. Recall the matrix description from
Proposition \ref{prop:primquot}
with matrix entries the bimodules $B_\subtorh^{\chi_i,\chi_j}\cong \ol{e_i}B\ol{e_j}$ for pairwise orthogonal idempotents
$\ol{e_i}$, $\ol{e_j}$.
This gives at least the $p$ column ideals $B\ol{e_i}$ in $B$. Since $B=B\ol{e_1}\oplus\ldots\oplus B\ol{e_p}$ is a direct sum
decomposition of $B$ it remains to check that inside the column ideals there are no further direct sums of left ideals
$I_i\oplus I_i\p$ of $B$. Write $B_{ij}$ for $B_\subtorh^{\chi_i,\chi_j}$.

Assume $B\ol{e_i}\supseteq I_i\oplus I_i\p$, we have to show that one of them must have been trivial. Assume $I_i$ is
    nontrivial.
We can find a regular element in $I_i\oplus\bigoplus_{j\neq i} B\ol{e_j}$:

First construct an element of $I_i$ that has exactly one nonzero entry located on the `diagonal' (recall that $I_i$ is
    contained in the $i^{\text{th}}$ column ideal). Let $a\in I_i$ be any nonzero element.
If the `diagonal entry' $a_{ii}$ is nonzero, multiply $a$ from the left with the idempotent $\ol{e_i}$, this gives the
desired element. In the case that $a_{ii}=0$, some other entry $a_{ji}$ must have been nonzero. By multiplication with
$\ol{e_j}$ we may assume that this was the only nonzero entry of $a$. We want to see $B_{ij}a_{ji}\neq0$, in other words,
one can find $b\in B_{ij}$ such that $c=ba \in B_{ii}$ has a nonzero `diagonal entry' at position $(i,i)$.
\begin{lemma}
If $B_{ij}a_{ji}=0$, then $a_{ji}$ was zero.
\end{lemma}
\begin{bew}
From $B_{ij}a_{ji}=0$ it follows immediately that $(a_{ji}B_{ij})(a_{ji}B_{ij})=0$ in the domain $B_{jj}$, hence
$a_{ji}B_{ij}=0$. Now consider the idempotent $e=\ol{e_i}+\ol{e_j}$. The matrices in $eBe$ only have nonzero entries at the
positions $(i,i)$, $(i,j)$, $(j,i)$, $(j,j)$. One calculates for $a$ as above the identities $ ((eBe)a(eBe))^2 =0$ and
	$$B(eBe)a(eBe)B\cdot B(eBe)a(eBe)B\ =\ B((eBe)a(eBe))^2 B\ = 0.$$
Since $B$ is a prime ring, $B(eBe)a(eBe)B=0$ follows. Moreover $a= e^4\cdot a\cdot e^4\in B(eBe)a(eBe)B=0$, so $a=0$ and in
particular $a_{ji}=0$.
\end{bew}
These arguments show that there is $c\in I_i$ with the $i^{\text{th}}$ entry $\neq 0$ and all other entries $=0$. One can
complete $c$ to a regular element in $I_i\oplus\bigoplus_{j\neq i} B\ol{e_j}$ by defining a matrix $C$ with $c_{ii}=c$, all
other diagonal entries are arbitrary regular elements from the domains $B_{jj}$ and all non-diagonal entries are zero. A
short calculation which can be found in \cite[Proof of Theorem 7.5.1]{dipl} shows that such $C$ is indeed regular in $B$.

From the existence of a regular element in the ideal $I_i\oplus\bigoplus_{j\neq i} B\ol{e_j}$ it follows that the
    ideal is essential in $B$, cf. \cite[Theorem 11.13, (2)$\Rightarrow$(5)]{LamLectures}. In particular, the intersection
    of $I_i\oplus\bigoplus_{j\neq i} B\ol{e_j}$ with $I_i\p$ is nonzero if $I_i\p$ is nonzero. But $I_i\oplus I_i\p$ was a
    direct sum, so $I_i\p=0$ follows.

Altogether we obtain
$$\grg{}\left(B_\subtor^\chi/J(\alpha)\right) =p= \text{ number of connected components of }\ol{\brac{\alpha}}_\bchi,$$
where we employed for the last equality Remark \ref{bem:zshgkomp}. The theorem follows.

\subsection{Ehrhart theory applied to families of regions $\ol{\brac{\alpha}}_\bchi$}\label{ssec:ehrhartalpha}

	We have seen in Remark \ref{bem:zshgkomp} that
	$\ol{\brac{\alpha}}_\bchi=\bigcup E\p+\delta_i+\alpha=\bigcup V(\subtorh-\chi_i(\subtorh))$.
	Here, $E\p$ is connected in the Zariski topology, and the translates do not intersect each other.
	Therefore, the connected components agree with $E\p+\delta_i+\alpha$, that is, with $V(\subtorh-\chi_i(\subtorh))$.
	Eventually one may calculate the Goldie rank of $B_\subtor^\chi/J(\alpha)$ via counting the translates.
	We reduce this task now to counting lattice points in an appropriate polytope, which can be done by Ehrhart
quasi-polynomials.

\begin{defi}
A {\it polyhedron} $P\subset\RR^n$ is given as the intersection of $m$ half-spaces,
$$P=P(A,z)= \{ x\in \RR^n\ |\ Ax\leq z, \text{ i.e. } (Ax)_i \leq z_i \ \forall 1\leq i \leq m \}.$$
A {\it polytope} is a bounded polyhedron. Its dimension is the dimension of the subspace of $\RR^n$ it spans. The {\it dilation
of a polytope} $P$ by $x$ is defined as $ xP\ =\ \{ x\cdot p\in\RR^n\ |\ p\in P\}$. A (full) lattice $L$ inside $\RR^n$ is given
by the $\ZZ$-span of some basis of $\RR^n$, the standard lattice is the one generated by the standard basis. For further details
see \cite{Ziegler} and \cite{Beck}.
\end{defi}

Classical Ehrhart theory states that given a lattice $L$ and an integral polytope $P$ of dimension $n$ (i.e. one with the
vertices on the lattice), the number of lattice points inside dilations $xP$ of $P$ by $x\in\ZZ_{>0}$ is a polynomial of degree
$n$ in $\QQ[x]$, see \cite[Theorem 3.8]{Beck} or \cite{ehrhart}.
This so-called {\it Ehrhart polynomial} is denoted by $\EHP_P(x)$.

The classical concept may be varied for rational polytopes (with vertices on the $\QQ$-span of the lattice). The lattice point
enumerator is then a quasi-polynomial in the dilation factor $x$, also denoted by $\EHP_P(x)$ \cite[Theorem 3.23]{Beck}. A {\it
quasi-polynomial} of degree $n$ and period $m$ is a polynomial of degree $n$ with coefficients that are periodic functions with
periodicity $m$. Since we are only interested in evaluating $\EHP_P$ for integral $x$, we may think of a quasi-polynomial as a
finite family of proper polynomials.

How does this apply to our situation?
Instead of considering $xP$ as a dilation of $P$, one can go one dimension up and interpret $xP$ as the intersection of a cone
with parallel translates of a hyperplane.
This should remind one of the intersection of the cone $\Delta_\alpha$ with $V(\subtor-\chi(\subtor))$, although we still have
to make some effort to tighten the analogy.
\begin{bem}
Notice that in Ehrhart theory one works over $\RR$, but due to our rationality assumptions on the arrangement this doesn't cause
any problem: We can count the lattice points inside $\QQ^n \subset k^n$ without losing information (and induce from $\QQ^n$ to
$\RR^n$ to apply Ehrhart theory).
\end{bem}
So, where do such families of integral or rational polytopes $P$ together with their dilates $xP$ for $x\in\ZZ_{>0}$ occur?
We encounter them when counting the connected components of $\ol{\brac{\alpha}}_\bchi$ for all $\alpha\in\Torus\ua$ with the
same sign configuration $J=J_+\cup J_-$.
Recall from Section \ref{ssec:closurereg} the description of $\ol{\brac{\alpha}}_\bchi$ as the intersection
$\ol{\brac{\alpha}}_\bchi=M_J\cap M_\vartheta$ where
\begin{align*}
M_J\ &=\ \left\{ \gamma\in V(\subtor-\chi(\subtor))\ \left|\
		\begin{matrix}
			\text{for all }j\in J\text{ one has }\gamma_j\in\ZZ,\\
			\gamma_j\geq0\text{ for }j\in J_+,\\
			\gamma_j<0\text{ for }j\in J_-
		\end{matrix}\right.\right\}\\
M_\vartheta\ &=\	\left\{ \gamma\in V(\subtor-\chi(\subtor))\ \left|\
		\begin{matrix}
			\sum\limits_{j\in J}\gamma_j\eta_j\ \in\ \sum\limits_{j\in J}\alpha_j\eta_j\ +\ \sum\limits_{i\in I}\ZZ\eta_i
		\end{matrix}\right.\right\}.
\end{align*}

\begin{bem}
	The set $M_J$ is highly appropriate for a description via lattice points and polytopes,
	whereas the behavior of $M_\vartheta$ can hardly be controlled.
	The problem is the missing description of $M_\vartheta$ in terms of lattices: We do not know how
	$M_\vartheta$ changes when changing $\alpha$.
\end{bem}

Since $M_\vartheta$ cannot be controlled, the following considerations concentrate on a handy special case. From now on, we make
the following assumptions:
\begin{enumerate}[label=\bf{Ass\arabic{*}}, ref=\bf{Ass\arabic{*}}]
\item\label{rat} Fix as usual $\subtor$ and thus $V(\subtor)$ \emph{rational}, i.e. the coefficients of the defining
    equations are in $\QQ$.
\item Fix a sign configuration $J=J_+\cup J_-$ for $\subtor$.
\item\label{mj} Assume
$\subtor\ua=\spann{k}\{\eta_j\ |\ j\in J\}\oplus \spann{k}\{\eta_i\ |\ i\in I\}.$
This guarantees $M_\vartheta=V(\subtor-\chi(\subtor))$, as can be easily checked, and therefore
$\ol{\brac{\alpha}}_\bchi=M_J$.
\end{enumerate}
Notice that we only fix the sign configuration $J$, not $\alpha$ itself, nor even $\chi$. We call $\alpha$ satisfying \eqref{rat}-\eqref{mj} of \emph{type J} if the associated sign configuration is $J$. We want to compare
$\ol{\brac{\alpha}}_\bchi$ with $\ol{\brac{\alpha\p}}_\bchip$, where $\alpha\p=x\cdot\alpha$ is a dilation of $\alpha$ with
$x\in\ZZ_{>0}$ and $\chi\p=\sum_{k=1}^n\alpha_k\p\eta_k$. This makes sense as long as $x\alpha$ belongs to the same index set $J$: If some of the coordinates of $\alpha$ are rational, then $x\alpha$ might have more integral coordinates than $\alpha$. In this case there are more inequalities activated in the definition of $\brac{x\alpha}_{B^{x\chi}}$, so the index set $J$ could also change. This we exclude by assuming that $x\alpha$ belongs to the same index set $J$. The property $M_\vartheta=V(\subtor-\chi(\subtor))$ is conserved under dilation of $\alpha$ and $\chi$:
\begin{lemma}
	Consider $\alpha\in\Torus\ua$ with $M_\vartheta=V(\subtor-\chi(\subtor))$
	and its dilation $\alpha\p=x\cdot\alpha$ ($x\in\ZZ_{> 0}$).
	Then one also has $M_\vartheta\p=V(\subtor-\chi\p(\subtor))$ for $\chi\p=\sum\limits_{k=1}^n\alpha_k\p\eta_k$.
\end{lemma}
From assumption \eqref{mj} it follows immediately that
$$\ol{\brac{\alpha}}_\bchi=M_J=\left\{ \gamma\in V(\subtor-\chi(\subtor))\ \left|\
		\begin{matrix}
			\text{for all }j\in J\text{ one has }\gamma_j\in\ZZ,\\
			\gamma_j\geq0\text{ for }j\in J_+,\\
			\gamma_j<0\text{ for }j\in J_-
		\end{matrix}\right.\right\}.$$
Counting the connected components resp. the different $\delta_k$ in
\begin{align}\label{eq:delta}
\ol{\brac{\alpha}}_\bchi=\bigcup_{k=1}^p E\p+(\delta_k+\alpha),
\end{align}
is equivalent to counting the elements in
$$D=\left\{ \delta\in V(\subtor-\chi(\subtor))\ \left|\
		\begin{matrix}
		\text{for all }i\in I\text{ one has }\delta_i=\alpha_i,\\
			\text{for all }j\in J\text{ one has }\delta_j\in\ZZ,\\
			\delta_j\geq 0\text{ for }j\in J_+,\\
			\delta_j<0\text{ for }j\in J_-
		\end{matrix}\right.\right\},$$
as the following lemma says:
\begin{lemma}\label{lem:deltadelta}
\begin{enumerate}[label=\roman{*}), ref=(\ref{lem:deltadelta}.\roman{*})]
\item $\delta+E\p\subset \ol{\brac{\alpha}}_\bchi$ for all $\delta\in D$ (hence $\delta$ represents one of the translates).
\item Every translate of $E\p$ in $\ol{\brac{\alpha}}_\bchi$ has a representative in $D$.
\item Let $\delta,\delta\p\in D$. Then $E\p+\delta=E\p+\delta\p$ if and only if $\delta=\delta\p$.
\end{enumerate}
\end{lemma}
\begin{bew}
\begin{enumerate}[label=\roman{*})]
\item is for free.
\item Consider an arbitrary element $\gamma\in \ol{\brac{\alpha}}_\bchi$. We can add any element $e\in E\p$ to $\gamma$
    without leaving $\ol{\brac{\alpha}}_\bchi$ or even the component where $\gamma$ lives. Define $e$ by $e_j=0$ for all
    $j\in J$, $e_i=\alpha_i-\gamma_i$ for all $i\in I$. This is in $E\p$ and $\gamma+ e\in D$ (use $\ol{\brac{\alpha}}_\bchi
    =M_J$).
\item The elements $\gamma$ in $E\p+\delta$ look like $e+\delta$ for some $e\in E\p$. From the definition of $E\p$ we know
    that $e_j=0$ and hence $\gamma_j=\delta_j$ for all $j\in J$. If $E\p+\delta=E\p+\delta\p$ we deduce
    $\delta_j=\delta_j\p$ for all $j\in J$. From the definition of $\delta,\delta\p\in D$ we have
    $\delta_i=\delta_i\p=\alpha_i$ for all $i\in I$, so $\delta=\delta\p$.
\end{enumerate}
\end{bew}
\begin{cor}
The number of connected components of $\ol{\brac{\alpha}}_\bchi$ equals the cardinality of $D$.
\end{cor}
Now that we have a set of representatives $D$ with fixed coordinates $\delta_i=\alpha_i$ for all $i\in I$, we may project from
$k^n$ onto $k^J=\spann{k}\{\pi_j\ |\ j\in J\}$ without changing the cardinality of $D$. This projection gives us the final pair
of a lattice and a polytope to which we can apply Ehrhart theory:
\begin{lemma}\label{lem:polyalpha}
Denote the projection $k^n\srar k^J$ by $\text{pr}_J$.
Then
\begin{enumerate}[label=\roman{*}), ref=(\ref{lem:polyalpha}.\roman{*})]
\item $\text{pr}_J(D)\ =\ \text{pr}_J\left(\ol{\brac{\alpha}}_\bchi\right) = P_J\cap L$, where $L\ =\ \ZZ^J\ =\
    \spann{\ZZ}\{\pi_j\ |\ j\in J\}$ is the standard lattice in $\QQ^J=\spann{\QQ}\{\pi_j\ |\ j\in J\}$ and
$P_J$ is the polyhedron defined by
\begin{align*}
P_J&=\left\{ d\in V_J\ \left|\
		\begin{matrix}
			d_j\geq 0&\text{ for }j\in J_+,\\
			d_j\leq-1&\text{ for }j\in J_-
		\end{matrix}\right.\right\}\\
\text{where }V_J&=\left\{v\in \QQ^J\ |\ \sum\limits_{j\in J}v_j\eta_j =\chi_J\right\}\ \text{ with }
\chi_J= \sum\limits_{j\in J}\alpha_j\eta_j\ \text{ and }\
\alpha_J= (\alpha_j)_{j\in J}\in\QQ^J .
\end{align*}
\item $P_J=\Delta_J\cap V_J$ has a description as the intersection of the cone
\begin{align*}
	\Delta_J\ &=\ \{v\in \QQ^J\ |\ v_j\geq 0\text{ for }j\in J_+,\ v_j\leq-1\text{ for }j\in J_-\}\\
	&=\ \left\{ v\in \QQ^J\ \left|\
		\begin{matrix}
			\Lambda_j(v)\leq 0\text{ for }j\in J_+,\\
			\Lambda_j(v)\leq -1\text{ for }j\in J_-
		\end{matrix}\right.\right\}\\
	\Lambda_j\ &=\
	\begin{cases}\Lambda_j(v)= -v_j,\ &\text{if }j\in J_+\\
	\Lambda_j(v)=v_j,\ &\text{if }j\in J_-\end{cases}
	\end{align*}
	with the affine subspace $V_J\subset\QQ^J$.
\item $P_J$ is a polytope, i.e. $P_J$ is a \emph{bounded} polyhedron.
\end{enumerate}
\end{lemma}
\begin{bew}
\begin{enumerate}[label=\roman{*})]
\item $\text{pr}_J(D) = P_J\cap L$ follows from keeping track of the definitions.
\item $P_J=\Delta_J\cap V_J$ follows also from the definitions.
\item $P_J$ is indeed a polyhedron since it has a description via linear inequalities (given in the definition of
    $\Delta_J$) and linear equalities (given in the definition of $V_J$). Now let us check the boundedness by showing that
    the translated polyhedron $P_J-\alpha_J = (\Delta_J-\alpha_J)\cap V_0$ is bounded. Here $V_0=\{v\in\QQ^J\ |\
    \sum\limits_{j\in J} v_j\eta_j=0\}$ is the image of the projection $\text{pr}_J(V(\subtor)\cap\QQ^n)$.
Notice that $\text{pr}_J(V(\subtor)\cap \QQ^n)= V_0$ and $\lambda_j(\gamma)=\Lambda_j(\text{pr}_J(\gamma))$ for all
$\gamma\in  V(\subtor)\cap \QQ^n$. Thus the image of $\Lambda_j$ on $V_0$ is determined by $\lambda_j$.
By Proposition \ref{lem:techgeo} and the definition of $J$, for $\lambda_j\in (V(\subtor)\cap\QQ^n)\ua$ there exist rational
coefficients $z_j>0$ such that
$$ \sum\limits_{j\in J} z_j\lambda_j =0.$$
For $\Lambda_j$ restricted to $V_0$ considered as an element of $V_0\ua$, this implies
$$ \sum\limits_{j\in J} z_j\Lambda_j =0.$$
Using the defining inequalities of $\Delta_J-\alpha_J$, we find upper and lower bounds for all coefficients $v_j$ of
$v\in(\Delta_J-\alpha_J)\cap V_0$:
\begin{itemize}
	\item $j\in J_+$: In this case $v_j\geq -\alpha_j$ and
	$$v_j =-\Lambda_j(v) = \frac{1}{z_j}\sum\limits_{k\neq j} z_k\Lambda_k(v) \leq \frac{1}{z_j}\sum\limits_{k\neq j}
z_k(-\alpha_k). $$
	\item $j\in J_-$: In this case $v_j\leq -1-\alpha_j$ and
	$v_j =\Lambda_j(v)= - \frac{1}{z_j}\sum\limits_{k\neq j} z_k\Lambda_k(v)\geq \frac{1}{z_j}\sum\limits_{k\neq j}
z_k\alpha_k.$
\end{itemize}
\end{enumerate}
\end{bew}
Notice that the definition of the polyhedral cone $\Delta_J$ is independent of $\chi$, the $\chi$ influences only $V_J$.
Now we have nearly reached our goal: We have
\begin{align*}
\grg{}(\bchi/J(\alpha))\ &=\ \text{number of connected components of }\ol{\brac{\alpha}}_\bchi\\
&=\ \text{cardinality of }\text{pr}_J\left(\ol{\brac{\alpha}}_\bchi\right)\\
&=\ \text{cardinality of }P_J\cap L
\end{align*}
where $P_J$ is a polytope and $L$ is a lattice, which enables us to apply Ehrhart theory. The only missing step is to describe
the set $\text{pr}_J\left(\ol{\brac{\alpha\p}}_\bchip\right)$ as `dilation' of
$\text{pr}_J\left(\ol{\brac{\alpha}}_\bchi\right)$, where $\alpha\p=x\alpha$ with $x\in\ZZ_{>0}$. Unfortunately, it is
\emph{not} the case that $\text{pr}_J\left(\ol{\brac{\alpha\p}}_\bchip\right)=xP_J\cap L$, since the center of dilation is not
necessarily the origin.
In fact, it is the point $z\in\QQ^J$ with coordinates
\begin{align}\label{eq:z}
z_j=\begin{cases}0\ &\text{for }j\in J_+\\ -1\ &\text{for }j\in J_-\end{cases}
\end{align}
depending on $J=J_+\cup J_-$ but not on $\chi$. In other words, one has to translate $P_J\cap L$ by $-z$, then dilate by an
appropriate dilation factor and translate the result back by $z$ to realize
$\text{pr}_J\left(\ol{\brac{\alpha\p}}_\bchip\right)$ from $\text{pr}_J\left(\ol{\brac{\alpha}}_\bchi\right)$.
The next proposition tells us how $x$ could be rescaled. It is rather technical, but notice that the rescaling is linear in $x$.
\begin{prop}[The dilation factor]\label{prop:dila}
	Let $\alpha\neq z$; so the polytope $P_J$ consists of more than one point.
	The description of $\text{pr}_J\left(\ol{\brac{\alpha}}_\bchi\right)$ as intersection $P_J\cap L$
	is up to a linear rescaling compatible with dilation, i.e. for $x\in\ZZ_{>0}$ such that $x\alpha$ induces the same configuration $J$ of signs,
	$$\text{pr}_J\left(\ol{\brac{x\alpha}}_\bchip\right) \ =\ \left(f(x)(P_J-z)+z\right)\cap L.$$
	Here we use the notation of Lemma \ref{lem:polyalpha},
	$$f(x)=\frac{x-a_0}{1-a_0}\qquad\text{where }a_0\in\QQ\text{ such that }z\in a_0\cdot\alpha_J +V_0$$
	and $V_0 = \{v\in\QQ^J\ |\ \sum\limits_{j\in J}v_j\eta_j =0\}$ (the analogue of $V(\subtor)$).
	The factor $a_0$ is rational and depends only on $\alpha$ and $z$, not on $x$.
\end{prop}
We have $a_0\neq 1$ because $\alpha\neq z$.
\begin{bew} The proof is a straightforward calculation.
	Notice that $\alpha\p=x\alpha$ and $\alpha$ are by assumption of the same type. Thus,
	$$\ol{\brac{\alpha\p}}_\bchip = M_J\p = \left\{ \gamma\in V(\subtor-\chi\p(\subtor))\ \left|\
		\begin{matrix}
			\text{for all }j\in J\text{ one has }\gamma_j\in\ZZ,\\
			\gamma_j\geq0\text{ for }j\in J_+,\\
			\gamma_j<0\text{ for }j\in J_-
		\end{matrix}\right.\right\}.$$
	Furthermore, by Lemma \ref{lem:polyalpha} we have $\text{pr}_J\left(\ol{\brac{\alpha\p}}_\bchip\right) =P_J\p\cap L$
	where
	\begin{eqnarray*}
	P_J\p=\left\{ d\in V_J\p\;\left|\;
		\begin{matrix}
			d_j\geq 0&\text{ for }j\in J_+,\\
			d_j\leq-1&\text{ for }j\in J_-
		\end{matrix}\right.\right\}&\text{and}&
	V_J\p =\left\{v\in \QQ^J\;|\;\sum\limits_{j\in J}v_j\eta_j =\sum\limits_{j\in J}\alpha\p_j\eta_j=\chi_J\p\right\}.
	\end{eqnarray*}
	We have to show that $P_J\p=f(x)(P_J-z)+z$ or equivalently
	$$P_J\p-z=f(x)(P_J-z).$$
	First, we determine $f(x)$.
	It should satisfy that multiplication of $(\alpha_J-z)$ by $f(x)$
	with subsequent translation by $z$ is equal to the original $x\alpha_J$ up to addition of some
	$v\in V_0= \left\{v\in\QQ^J\ |\ \sum_{j\in J}v_j\eta_j =0\right\}.$
	In other words, we want $f(x)(\alpha_J-z)+z=x\alpha_J +v.$ Take $a_0\in\QQ$ such that
	$z\in a_0\cdot\alpha_J\in V_0$, i.e. $z=a_0 \alpha_J + v_0$ with $v_0\in V_0$.
	A short calculation gives
	$$f(x)= \frac{x-a_0}{1-a_0} \quad\text{and}\quad v= \frac{1-x}{1-a_0}v_0.$$
	This does the job: One can calculate that
	\begin{align*}
	f(x)(P_J-z)\ &=\ \left\{ d\in \QQ^J\ \left|\
		\begin{matrix}
			\sum\limits_{j\in J}d_j\eta_j &= f(x)\cdot \sum\limits_{j\in J}(\alpha_j-z_j)\eta_j\\
			d_j\geq 0&\ \text{ for }j\in J_+\\
			d_j\leq 0&\ \text{ for }j\in J_-
		\end{matrix}\right.\right\}.
	\end{align*}
	and
	\begin{align*}
	P_J\p-z\ &=\ \left\{ d\in \QQ^J\ \left|\
		\begin{matrix}
			\sum\limits_{j\in J}d_j\eta_j &= \sum\limits_{j\in J}(x\alpha_j-z_j)\eta_j\\
			d_j\geq 0&\ \text{ for }j\in J_+\\
			d_j\leq 0&\ \text{ for }j\in J_-
		\end{matrix}\right.\right\}.
	\end{align*}
	Comparing them it remains to check whether
	$$f(x)\cdot \sum\limits_{j\in J}(\alpha_j-z_j)\eta_j=\sum\limits_{j\in J}(x\alpha_j-z_j)\eta_j.$$
	But this holds since $f(x)(\alpha_J-z)=(x\alpha_J-z)+v$ and $v\in V_0$, i.e. $\sum_{j\in J} v_j\eta_j=0$.
\end{bew}
Luckily, $z$ is itself a lattice point, so that we do not need to translate the arrangement back by $z$ to count the lattice
points. The benefit of the whole construction is the following theorem:
\begin{thm}\label{thm:zshgzahl}
For \emph{all} regions of the form $\ol{\brac{x\alpha}}_{B^{x\chi}}$ with $x\in\ZZ_{>0}$ such that $x\alpha$ has the same type $J$, the number of connected components
equals $\#(f(x)(P_J-z)\cap L)$, where $f(x)$ is the rescaling of the dilation factor of Proposition \ref{prop:dila}.
\end{thm}

\subsection{Goldie ranks using Ehrhart quasi-polynomials}\label{ssec:goldpol}

Recall the assumptions Ass1-Ass3 from the previous section.
Now choose $\alpha\in\Torus\ua$ of type $J=J_+\cup J_-$ (in particular, $\alpha_j\in\ZZ$ for all
$j\in J$).
We concluded the last section by the insight that the Goldie rank of the primitive quotient $\bxchi/J(x\alpha)$ with
$x\in\ZZ_{>0}$ such that $x\alpha$ is of also of type $J$ may be expressed as the number of lattice points $\#(f(x)(P_J-z)\cap L)$ where $L=\ZZ^J$ is the standard lattice
in the ambient space $\QQ^J$ and $P_J$ is an appropriate polytope, translated by $z$ (only depending on $J$, see \eqref{eq:z}).
It is important to notice that $P_J$ is a rational polytope since due to assumption \eqref{rat} the defining equations in the
definition of $P_J$ have rational coefficients, see Lemma \ref{lem:polyalpha}. The dilation factor $f(x)\in\QQ$ is a linear
rescaling of $x$, only dependent of the fixed $\alpha$ and $z$.

Now one would like to compute the Goldie rank via the Ehrhart quasi-polynomial, i.e.
$$\#(f(x)(P_J-z)\cap L)=\EHP_{P_J-z}(f(x)).$$
Unfortunately to apply the classical Ehrhart theory we need $f(x)\in\ZZ$. Alternatively we can work with Ehrhart
quasi-polynomials for \emph{rational} dilation factors as defined in \cite{linke}:
\begin{thm}[Goldie rank of $\bchi/J(\alpha)$, first version]
\label{thm:goldpol1}
Let $\alpha\in\Torus\ua$ of type $J=J_+\cup J_-$ be given. Let $x\in\ZZ_{>0}$ such that $x\alpha$ is also of type $J$. Assume conditions
\eqref{rat}-\eqref{mj} are satisfied. The Goldie rank of the primitive quotients $\bxchi/J(x\alpha)$ is then a quasi-polynomial
in $x$ given by
$$\grg{}(B^{x\chi}/J(x\alpha))\ =\ \EHP^\QQ_{P_J-z}(f(x)),$$
where $\EHP^\QQ_{P_J-z}(f(x))$ is the rational Ehrhart quasi-polynomial of the rational polytope $P_J-z$
with respect to the standard lattice $L=\ZZ^J$ composed with the linear rescaling $f(x)\in\mathbb{Q}$ from Proposition~\ref{prop:dila}.
\end{thm}
\begin{bew}
By Theorem \ref{thm:conncpt} and  Theorem \ref{thm:zshgzahl} we know that $\grg{}(B^{x\chi}/J(x\alpha))$ equals the number of
connected components of $\ol{\brac{x\alpha}}_{B^{x\chi}}$ which is equal to $\#(f(x)(P_J-z)\cap L)$. By \cite[Theorem
1.2]{linke}, this is given by the rational Ehrhart quasi-polynomial $\EHP^\QQ_{P_J-z}$ in $f(x)$.
\end{bew}
In particular, for fixed $\alpha$, we have a finite family of genuine polynomials that compute the Goldie ranks of
$B^{x\chi}/J(x\alpha)$.

However, we do not need to use the Ehrhart quasi-polynomial for arbitrary rational dilation factors. Instead we may modify the
polytope in such a way that precomposition with $f$ is encapsulated in the definition of the `better' polytope. In our special
situation, the dilation factor is of the form $\frac{x-a_0}{1-a_0}$ with $a_0=\frac{a_Z}{a_N}\in\QQ$. From the construction of
$a_0$ we read off that $a_0<1$, i.e. $a_N-a_Z>0$, and one gets
$$\frac{x-a_0}{1-a_0} = \frac{a_N x+ a_Z}{a_N-a_Z}.$$
Since $a_0$ was independent of $x$, so are $a_N$ and $a_Z$.
Define the new reference polytope $Q= \frac{1}{a_N-a_Z} (P_J-z)$.
Hence
$$f(x)\cdot (P_J-z)\ =\ (a_N x+ a_Z)\cdot Q,$$
converting the rational dilation of $P_J-z$ to integral dilation of $Q$.
This gives
\begin{thm}[Goldie rank of $\bchi/J(\alpha)$, second version]\label{thm:goldpol}
Let $\alpha\in\Torus\ua$ of type $J=J_+\cup J_-$ and $x\in\ZZ_{>0}$ such that $x\alpha$ is again of type $J$. Assume conditions
\eqref{rat}-\eqref{mj} are satisfied. Then the Goldie rank of the primitive quotients $B^{x\chi}/J(x\alpha)$ is a
quasi-polynomial in $x$ given by
$$\grg{}(B^{x\chi}/J(x\alpha))\ =\ \EHP_{Q}(a_N x+ a_Z),$$
where $\EHP_{Q}(a_N x+ a_Z)$ is the classical Ehrhart quasi-polynomial of the rational polytope $Q = \frac{1}{a_N-a_Z} (P_J-z)$
with respect to the standard lattice $L=\ZZ^J$ composed with the \emph{integral} linear rescaling $a_N x+ a_Z$.
\end{thm}

\begin{bem}
Let us briefly analyse the condition that $x\alpha$ is of the same type $J$:
It assures that the inequalities for the description of $\ol{\brac{x\alpha}}_{B^{x\chi}}$ are the same as for $\ol{\brac{\alpha}}_\bchi$. Therefore, we work inside the same cone $\Delta_J$. This means that we only admit those $x$ that are not divided by the denominators of the rational coordinates of $\alpha$. In particular, if $\alpha$ is integral, then $x\alpha$ is of the same type and satisfies \eqref{rat}-\eqref{mj} automatically. On the other hand, using this assumption, we can drop the requirement $x\in\ZZ_{>0}$ and instead formulate the result for any $x\in\QQ$ as long as the sign configuration remains the same for $x\alpha$.
\end{bem}
\begin{bem}
Note that in the construction more than one different quasi-polynomial occurs, but only finitely many altogether. Notice that we have one quasi-polynomial per cone $\Delta_J$ (of which there are $2^{|J|}$) and at most $2^n$ different index sets $J$. For given $\alpha$, which of the quasi-polynomials occur can be expressed by the denominators of the rational coordinates of $\alpha$ that divide $x$. There is only one caveat: Assumptions \eqref{rat}-\eqref{mj} must be satisfied for $x\alpha$.
\end{bem}

\begin{bem}
A result similar, but in fact easier, to Theorem \ref{thm:goldpol} holds for the primitive quotients of a central quotient of the hypertoric enveloping
algebra \cite[Proposition 7.4, Remark 7.5]{blpw}.
\end{bem}
Using well-known results about Ehrhart quasi-polynomials (see e.g. \cite{Beck}), one may deduce further properties of the Goldie
rank quasi-polynomials:
\begin{enumerate}[label=\roman{*})]
\item The degree of the Goldie rank quasi-polynomial is at most $|J|$, see \cite[Theorem 3.23]{Beck}.
\item The period of the quasi-polynomial, i.e. the number of genuine polynomials forming the Goldie rank quasi-polynomial,
    always divides
$\min\{d\in\ZZ_{>0}\ |\ dQ\text{ is an integral polytope}\}$, see \cite[Section 2.7, Theorem 3.23]{Beck}.
\item There are polynomial-time algorithms for the computation of the coefficients of an Ehrhart quasi-polynomial (in fixed
    dimension) available \cite[Chapter 5]{Barvinok}.
\end{enumerate}

\bibliographystyle{amsalpha}
\addcontentsline{toc}{section}{References}
\bibliography{literatur}

\end{document}